\documentclass[11pt,a4paper]{amsart}
\usepackage[latin1]{inputenc}
\usepackage[arrow, matrix, curve]{xy}
\usepackage{amssymb}
\usepackage{rotating}
\usepackage{amsmath}
\usepackage{latexsym}
\usepackage{xypic}
\usepackage{epic}
\usepackage{eepic}
\usepackage{graphicx}
\usepackage{longtable}
\usepackage{epsfig}
\usepackage{xypic}

\newtheorem{rem}{Remark}[section]
\newtheorem{prop}[rem]{Proposition}
\newtheorem{lemma}[rem]{Lemma}

\newtheorem{defi}[rem]{Definition}
\newtheorem{theorem}[rem]{Theorem}

\newcommand{\la}{\leftarrow}

\newcommand{\bprf}{{\it Proof.~}}
\newcommand{\binf}{{\it In fact }}
\newcommand{\ra}{\rightarrow}

\newcommand{\erem}{\hfill $\square$}

\newcommand{\N }{ \mathbb{N}}
\newcommand{\C }{ \mathbb{C}}
\newcommand{\R}{ \mathbb{R}}
\newcommand{\Z}{\mathbb{Z}}
\newcommand{\Q}{\mathbb{Q}}

\input epsf
\input cyracc.def

\makeatletter
\def\blfootnote{\xdef\@thefnmark{}\@footnotetext}

\begin{document}
\title{Symplectic automorphisms on Kummer surfaces}
\author{Alice Garbagnati}
\address{Alice Garbagnati, Dipartimento di Matematica, Universit\`a di Milano,
  via Saldini 50, I-20133 Milano, Italia}

\email{garba@mat.unimi.it}

\begin{abstract}
Nikulin proved that the isometries induced on the second
cohomology group of a K3 surface $X$ by a finite abelian group $G$
of symplectic automorphisms are essentially unique. Moreover he
computed the discriminant of the sublattice of $H^2(X,\Z)$ which
is fixed by the isometries induced by $G$. However for certain
groups these discriminants are not the same of those found for
explicit examples. Here we describe Kummer surfaces for which this
phenomena happens and we explain the difference.
\end{abstract}

\maketitle


\section{Introduction}
An automorphism $\sigma$ of finite order on a K3 surface $X$ is
symplectic if the desingularization  of the quotient, denoted by
$\widetilde{X/\sigma}$, is again a K3 surface. The finite abelian
groups of automorphisms on a K3 surface were classified by Nikulin
\cite{Nikulin symplectic}. Later Mukai \cite{Mukai} classified all
the finite groups of symplectic automorphisms. The main result of
\cite{Nikulin symplectic} is that the group of isometries of the
second cohomology group induced by a finite abelian group of
symplectic automorphisms of the K3 surfaces is essentially unique.
So these isometries do not depend on the surface but only on the
group. If the group $G$ of symplectic automorphisms is generated
by an involution, Morrison \cite{morrison} found this isometry
explicitly. A description of these isometries for the others
finite abelian groups acting symplectically on a K3 surface is
given in \cite{symplectic prime} and \cite{symplectic not prime}.
Since the isometries induced by $G$ are essentially unique, the
lattice $H^2(X,\Z)^G$ (the sublattice of the $G$-invariant
elements in $H^2(X,\Z)$) depends on $G$ but not on the K3 surface
$X$, up to isometry. In particular the discriminant of
$H^2(X,\Z)^G$ depends only on $G$. In \cite{Nikulin symplectic}
this discriminant is found using the relation between the K3
surface $X$ and the K3 surface $Y$ which is the desingularization
of $X/G$. This is an interesting idea essentially because $Y$
seems to be easier to study then $H^2(X,\Z)^G$. In fact
$H^2(Y,\Z)$ contains a sublattice $M_G$ generated, at least over
$\Q$, by rational curves with simple intersection properties,
which is related to $H^2(X,\Z)^G$ but is less complicated to
analyze. The construction used by Nikulin to find the discriminant
of $H^2(X,\Z)^G$ is summarized here in
\textit{Section \ref{section: Nikulin's construction}}.\\
Another way to compute the discriminant of $H^2(X,\Z)^G$ is to
construct an explicit example of a K3 surface $X$ admitting $G$ as
group of symplectic automorphisms and to find the lattice
$H^2(X,\Z)^G$ explicitly (this method is the one adopted in
\cite{symplectic prime}, \cite{symplectic not prime}). The
discriminants given in \cite{Nikulin symplectic} are not always
equal to the ones computed on explicit examples. Here we present
such example. Moreover we explain where the method described in
\cite{Nikulin symplectic} fails. In fact Nikulin asserts that if
$G$ is a symplectic group of automorphisms on $X$, then a certain
map between a lattice related to $Y$ and a lattice related to $X$
would be a surjective map. In {\it Section \ref{section: spectral
sequence and...}} we will prove that if the group $G$ is a cyclic
group, that map is indeed surjective, but in {\it Section
\ref{section: group of automorphisms of order two.}} we will
exhibit some counterexamples in the case $G$ is not cyclic. In
these counterexamples $X$ is a Kummer surface. The translations by
2-torsion points on an Abelian surface $A$ induce symplectic
involutions on the associated Kummer surface $Km(A)$. So the group
$(\Z/2\Z)^4$ is a group of symplectic automorphisms on each Kummer
surface. For this example we will explicitly describe the map
introduced by Nikulin.\\

{\it This paper is an excerpt from my PhD thesis. I would like to
thank Bert van Geemen, my thesis supervisor, for his unrivalled
assistance and backing. I would also like to thank the Mathematics
Department of the University of Milan and, in particular, Prof. A.
Lanteri, for the support and encouragement he provided me.}

\section{Nikulin's construction}\label{section: Nikulin's
construction} Let $X$ be a K3 surface. The second cohomology group
of $X$ with the cup product is a lattice isometric to the lattice
(which does not depend on $X$) $\Lambda_{K3}:=U\oplus U\oplus
U\oplus E_8(-1)\oplus E_8(-1)$, where $U$ is the rank two lattice
with pairing $\left[\begin{array}{rr}0&1\\1&0\end{array}\right]$
and $E_8(-1)$ is the rank eight lattice associated to the Dynkin
diagram
$E_8$ with the pairing multiplied by $-1$ (cf.\ \cite{bpv}).\\
Let $\sigma$ be an automorphism of $X$. Then $\sigma^*$ is an
isometry of $H^2(X,\Z)$ which preserves the set of the effective
divisor and its $\C$-linear extension preserves the Hodge
decomposition of $H^2(X,\C)$.
\begin{defi}
An automorphism $\sigma\in Aut(X)$ is \textbf{symplectic} if and
only if $\sigma^*_{|H^{2,0}(X)}=Id_{|H^{2,0}(X)}$. A group of
automorphisms acts symplectically on $X$ if all the elements of
the group are symplectic automorphisms.
\end{defi}

\begin{rem}{\rm \cite[Theorem 3.1 b)]{Nikulin symplectic} An automorphism $g$ of $X$ is symplectic if and only if
$g^*_{T_X}=Id_{T_X}$.\erem}\end{rem} Let $G$ be a finite group
acting symplectically on $X$. The quotient of $X$ by $G$ is a
singular surface $X/G$. Let us call $Y=\widetilde{X/G}$ the
minimal resolution of this surface. Then we have the following
diagram
\begin{eqnarray*}
\begin{array}{lcr}
X&\dashrightarrow& \widetilde{X/G}=Y\\
&^{^\pi{^\searrow}}X/G^{^{^{\ \swarrow}\ \beta}}
\end{array}
\end{eqnarray*}
Here we recall one of the main results on the symplectic
automorphisms of K3 surfaces.
\begin{theorem} Let $X$ be a K3 surface and $G$ be a finite group of
symplectic automorphisms of $X$. Then $Y=\widetilde{X/G}$ is a K3
surface.\\
Viceversa, let $X$ be a K3 surface and $G$ be a finite group of
automorphisms of $X$. If $Y=\widetilde{X/G}$ is a K3 surface, then
$G$ is a group of symplectic automorphisms.
\end{theorem}

From now on we will assume that the group $G$ is a finite abelian
group.

\begin{defi}\label{defi: the map theta} Let $X'=X\backslash\{\ldots,x_i,\ldots\}$ where
$\{\ldots,x_i,\ldots\}$ is the set of points with a non-trivial
stabilizer under the action of $G$. Let $Y'$ be the surface
$Y'=X'/G$ and $\theta:X'\ra Y'$ be the quotient map.
\end{defi}
\begin{rem}{\rm Since the action of $G$ on $X'$ is fixed points free,
the surface $Y'$ is smooth, and of course $Y'=Y\backslash \{\cup_j
M_j\}$ where $M_j$ are the curves which arise from the resolution
of the singularities of $X/G$. Moreover the map $\theta:X'\ra Y'$
is an unramified $|G|$-sheeted cover.\erem}\end{rem}
\begin{rem}{\rm Let $\vartheta$ be the rational map associated to the quotient $X/G$. So we have:
\begin{eqnarray}\label{formula: vartheta}
\begin{array}{rllllr}
&\widetilde{X}&        &&\\
&\downarrow   &\searrow&&\\
\vartheta:&X&\dashrightarrow &Y&\\
&\cup&&\cup\\
\phantom{aaaaaaaaaaaaaaaaaaaaaaaaaa}\vartheta_{|X'}=\theta:&X'&\ra&Y'&\phantom{aaaaaaaaaaaaaaaaaaaaaaaaaa}\square
\end{array}
\end{eqnarray}}
\end{rem}

\begin{defi}{\rm \cite[Definition 4.6]{Nikulin symplectic}}
We say that a group $G$ \textbf{has a unique action on
$\Lambda_{K3}$} if, given any two embeddings $i:G\hookrightarrow
Aut(X)$ and $i':G\hookrightarrow Aut(X')$ under which $G$ is a
group of symplectic automorphisms of the K3 surfaces $X$ and $X'$,
there exists an isometry $\varphi:H^2(X,\Z)\ra H^2(X',\Z)$ such
that $i'(g)^*=\varphi\circ i(g)^*\circ \varphi^{-1}$ for any $g\in
G$.
\end{defi}
\begin{theorem}\label{theorem: uniquess isometry, Nikulin 4.7}{\rm \cite[Theorem 4.7]{Nikulin symplectic}} Any
finite abelian group has a unique action on $\Lambda_{K3}$.
\end{theorem}
The natural question posed by this result is to find explicitly
these actions on $\Lambda_{K3}$. If the group acting
symplectically is $G=\Z/2\Z$, i.e. it is generated by a symplectic
involution, the answer to this question was found by Morrison in
\cite{morrison}. The cases $G=\Z/p\Z$ where $p$ is a prime are
analyzed in \cite{symplectic prime} and the other cases are
analized in \cite{symplectic not prime}. To find these isometries
the main problem is to find the sublattices $H^2(X,\Z)^G$ and
$(H^2(X,\Z)^G)^{\perp}$. Nikulin does not determine them, but
describes a technique to compute some invariants of these
lattices, which we recall here.

\begin{lemma}\label{lemma: H0, H1 of X,X',Y,Y'}
The surfaces $X$, $X'$, $Y$ and $Y'$ have the following
properties:\\ $\bullet$
$H^0(X,\Z)=H^0(X',\Z)=H^0(Y,\Z)=H^0(Y',\Z)=H_0(X,\Z)=H_0(X',\Z)=H_0(Y,\Z)=H_0(Y',\Z)=\Z$;\\
$\bullet$ $\pi_0(X')$ is trivial;\\
$\bullet$ $H_1(X,\Z)=H_1(X',\Z)=H_1(Y,\Z)=0$ and $H_1(Y',\Z)=G$;\\
$\bullet$ $H^1(X,\Z)=H^1(X',\Z)=H^1(Y,\Z)=H^1(Y',\Z)=0$.
\end{lemma}
\bprf All the surfaces $X$, $X'$, $Y$, $Y'$ are clearly
path connected and so their homology and cohomology group in degree 0 is $\Z$ and $\pi_0(X')$ is trivial.\\
The surfaces $X$ and $Y$ are K3 surfaces, so by definition
$H^1(X,\Z)=H^1(Y,\Z)=0$. The complex surface $X'$ is the surface
$X$ without some points. The four dimensional topological variety
$X$ is simply connected, this implies that $X'$ is simply
connected, i.e. $\pi_1(X')=0$. So $H_1(X',\Z)=Ab(\pi_1(X'))=0$.
Since $H^n(X',\Z)\simeq Hom(H_{n}(X'),\Z)\oplus
Ext(H_{n-1}(X'),\Z)$ 
and $Ext(\Z,\Z)=0$, we have
$H^1(X',\Z)=0$.\\
Since $X'$ is simply connected and $G$ acts without fixed points
on $X'$, it follows that $\pi_1(Y')\simeq G$.
Since $G$ is an abelian group, $H_1(Y',\Z)=G$. Hence
$H^1(Y',\Z)\simeq Hom(G,\Z)\oplus Ext (\Z,\Z)=0$.\erem
\begin{lemma}\label{lemma: H2Y' and H2X'} We have
$H^2(Y',\Z)\simeq H^2(Y,\Z)/\oplus_j \Z M_j$ and $H^2(X',\Z)\simeq
H^2(X,\Z)$ {\rm (see \cite[Lemma 6.1]{Nikulin symplectic},
\cite[Lemma 2]{Xiao}).}\end{lemma}

\begin{defi}\label{defi: MG}
Let $M_G$ be the minimal primitive sublattice of $\Lambda_{K3}$
containing the curves $M_j$ arising from the resolution of the
singularities of $X/G$. Let $P_G$ be its orthogonal with respect
to the bilinear form of $\Lambda_{K3}$.
\end{defi}
\begin{rem}\label{rem: MG is contained in NS(Y)}{\rm The lattice $M_G$ is primitively
embedded in $NS(Y)$, because it is generated over $\Q$ by the
curves $M_j$ on $Y$.\erem}\end{rem}
Since $M_G\oplus P_G\hookrightarrow \Lambda_{K3}$ with a finite
index and $\Lambda_{K3}$ is a unimodular lattice, the discriminant
group of $M_G$ is the same as the discriminant group of $P_G$ and
$M_G^{\vee}/M_G=P_G^{\vee}/P_G$, where $L^{\vee}=Hom_{\Z}(L,\Z)$
is the dual of the lattice $L$. Moreover we have the following
exact sequence
\begin{eqnarray*}
0\ra M_G\longrightarrow
\Lambda_{K3}\stackrel{l_b}{\longrightarrow}
P_G^{\vee}\longrightarrow 0
\end{eqnarray*}
where $b$ is the bilinear form on $\Lambda_{K3}$ and
$l_b(y)(x)=b(y,x)$, for $y\in \Lambda_{K3}$ and $x\in P_G.$ From
this sequence we obtain
\begin{eqnarray}\label{formula: Lambda/MG=PG}\Lambda_{K3}/M_G\simeq P_G^{\vee}.\end{eqnarray}
The map $\theta:X'\ra Y'$ induces a map $\theta^*:H^2(Y',\Z)\ra
H^2(X',\Z)^G$. Thanks to Lemma \ref{lemma: H2Y' and H2X'} we have
$\theta^*:H^2(Y,\Z)/\oplus \Z M_j=H^2(Y',\Z)\ra
H^2(X',\Z)^G=H^2(X,\Z)^G$. Let $\sigma_Y$ (resp. $\sigma_X$) be
the isometry between $H^2(Y,\Z)$ (resp. $H^2(X,\Z$)) and
$\Lambda_{K3}$. So we have the following exact sequences
\begin{equation}\label{diagram: sequences on H2X' and H2Y' }
\begin{array}{ccccccc}
0&\ra &tors(H^2(Y',\Z))&\ra &H^2(Y',\Z)&\xrightarrow{\ \ \theta^*\
\ }&
H^2(X',\Z)^G\\
& &\downarrow\!\! \wr& &\alpha\downarrow\!\! \wr&&
\beta\downarrow\!\! \wr\\
0&\ra &tors(H^2(Y,\Z)/\oplus\Z M_j)&\ra &H^2(Y,\Z)/\oplus \Z
M_j&\xrightarrow{\ \ \ \ \ \ }&
H^2(X,\Z)^G\\
& &\downarrow\!\! \wr& &\!\!\sigma_Y\!\downarrow\!\! \wr&&
\!\!\sigma_X\!\downarrow\!\! \wr\\
0&\ra &M_G/\oplus \Z M_j&\ra &\Lambda_{K3}/\oplus \Z
M_j&\xrightarrow{\ \ \ \ \ \ }& \Lambda_{K3}^G.
\end{array}
\end{equation}
From these sequences we obtain
\begin{eqnarray}\label{formula: L/M=Im}\Lambda_{K3}/M_G\simeq \mbox{Im}(\Lambda_{K3}/\oplus_j\Z M_j\ra \Lambda_{K3}^G)\simeq\theta^*(H^2(Y',\Z)).\end{eqnarray}
So $\Lambda_{K3}/M_G\simeq P_G^{\vee}$ is a submodule of
$\Lambda_{K3}^G$.\\
If $\theta^*$ were surjective, this would implies
$P_G^{\vee}\simeq\Lambda_{K3}/M_G\simeq \Lambda_{K3}^G$ (by
formulas \eqref{formula:
Lambda/MG=PG} and \eqref{formula: L/M=Im}).\\
We want to study the inclusion $P_G^{\vee}\subset \Lambda_{K3}^G$
and thus we need to study $\theta^*:H^2(Y',\Z)\ra H^2(X',\Z)$.
More precisely since $P_G^{\vee}\subset H^2(Y',\Z)$ and
$\Lambda_{K3}^G\simeq H^2(X,\Z)^G$ we need to consider
$\vartheta^*(P_G^{\vee})\subset H^2(X,\Z)^G$.
Since $P_G\subset H^2(Y',\Z)$,
$\vartheta^*(P_G^{\vee})=\theta^*(P_G^{\vee})$ (cf.\
\eqref{formula: vartheta}). Hence
\begin{equation}\label{formula: theta* is surjective iff}\theta^*\mbox{ is surjective if and only if
}H^2(X,\Z)^G\simeq \theta^*(P_G^{\vee}).\end{equation} As
lattices,  $\theta^*$ induces a scaling on
$P_G^{\vee}:\theta^*(P_G^{\vee})\simeq P_G^{\vee}(|G|)$ because
$(\theta^*(x),\theta^*(y))=|G|(x,y)$ for each $x,y\in (\oplus \Z
M_j)^{\perp}$ as $\theta:X'\ra Y'$ is of degree $|G|$.

\begin{rem}{\rm In \cite[Section 8,]{Nikulin symplectic} the author
states that $\theta^*$ is surjective. Assuming the surjectivity of
the map the discriminant of the group $\Lambda_{K3}^G$ is computed
in \cite[Lemma 10.2]{Nikulin symplectic} .\erem}\end{rem}

\section{Spectral sequences and the map $\theta^*:H^2(Y',\Z)\ra
H^2(X',\Z)^G$}\label{section: spectral sequence and...} We use a
spectral sequence to analyze the properties of the map $\theta^*$.
The main result is the sequence \eqref{formula: THE SEQUENCE}.
\begin{theorem}{\rm \cite[Theorem 6.10]{Weibel}}\label{theorem: convergence of E2pq to Hp+q(S/G,A)} Let $G$ be a group acting properly on
a space $S$ such that $\pi_0(S)$ is trivial. Then for every
abelian group $A$ there exists a spectral sequence and
$$E_2^{p,q}=H^p(G,H^q(S,A))\Rightarrow H^{p+q}(S/G,A).$$
\end{theorem}
\begin{rem}{\rm This sequence is a first quadrant
spectral sequence (i.e. $E_r^{p,q}=0$ if $p<0$ or $q<0$). Moreover
$E_2^{p,q}=0$ for each $q>\dim(S)$.\erem}
\end{rem}
Under the hypothesis of the theorem, the spectral sequence
converges, and so there are filtrations
$$E^2=E_0^2\supset E_1^2\supset E_2^2\supset E_3^2=0\ \ \mbox{ and
}\ \ E_{\infty}^{0,2}=E_0^2/E_1^2,$$ so we have a exact sequence $
0\ra E_1^2\ra E_0^2\ra E_{\infty}^{0,2}\ra 0.$ Analogously there
is an exact sequence $0\ra E_2^2\ra E_1^2\ra E_{\infty}^{1,1}\ra
0.$ So we have
\begin{eqnarray}\label{formula: diagrammone}
\begin{array}{ccccccccc}
&&0&&&&&&\\
&&\downarrow&&&&&&\\
&&E_2^2&&&&&&\\
&&\downarrow&&&&&&\\
0&\ra&E_1^2&\ra&E^2&\ra&E_{\infty}^{0,2}&\ra&0\\
&&\downarrow&&&&&&\\
&&E_{\infty}^{1,1}&&&&&&\\
&&\downarrow&&&&&&\\
&&0&&&&&&
\end{array}
\end{eqnarray}

We now use the same notation as in the previous section.\\
The hypothesis of Theorem \ref{theorem: convergence of E2pq to
Hp+q(S/G,A)} are satisfied by the topological space $S=X'$, the
group $G$ acting properly on $X'$ and $A=\Z$. We apply the Theorem
\ref{theorem: convergence of E2pq to Hp+q(S/G,A)} to the spectral
sequence $E_2^{p,q}=H^p(G,H^q(X',\Z))$. Since $H^1(X',\Z)=0$
(Lemma \ref{lemma: H0, H1 of X,X',Y,Y'}) we have
\begin{eqnarray}\label{formula: Ep1=0} E_2^{p,1}=H^p(G,H^1(X',\Z))=H^p(G,0)=0.\end{eqnarray}
We will compute the groups in \eqref{formula: diagrammone} in this
case.\\
\textbf{Computation of $E_{\infty}^{0,2}$.} Since the spectral
sequence is a first quadrant spectral sequence
$E_4^{0,2}=E_{\infty}^{0,2}$ ($E_r^{0,2}=E_{\infty}^{0,2}$ if
$r>\max\{p,q+1\}$). By definition
\begin{eqnarray}\label{formula: Einfty02 first} E_{\infty}^{0,2}=E_4^{0,2}=\frac{\ker(d_3^{0,2}:E_3^{0,2}\ra
E_3^{3,0})}{\mbox{Im}(d_3^{-3,4}:E_3^{-3,4}\ra
E_3^{0,2})}=\ker(d_3^{0,2}:E_3^{0,2}\ra E_3^{3,0})\end{eqnarray}
where $d_3^{-3,4}=0$ because the spectral sequence is in the first
quadrant. But $$E_3^{0,2}=\frac{\ker(d_2^{0,2}:E_2^{0,2}\ra
E_2^{2,1})}{\mbox{Im}(d_2^{-2,3}:E_2^{-2,3}\ra
E_2^{0,2})}=\ker(d_2^{0,2}:E_2^{0,2}\ra E_2^{2,1})=E_2^{0,2}$$
where the last equality is a consequence of \eqref{formula:
Ep1=0}. Moreover $$E_3^{3,0}=\frac{\ker(d_2^{3,0}:E_2^{3,0}\ra
E_2^{5,-1})}{\mbox{Im}(d_2^{1,1}:E_2^{1,1}\ra
E_2^{3,0})}=\frac{\ker(d_2^{3,0}:E_2^{3,0}\ra
0)}{{\mbox{Im}(d_2^{1,1}:0\ra E_2^{3,0})}}=E_2^{3,0}.$$
Substituting the last two formula in \eqref{formula: Einfty02
first} we obtain
$$
E_{\infty}^{0,2}=\ker(E_2^{0,2}\ra E_2^{3,0}).$$ By the definition
of the spectral sequence and by the equality
$H^0(G,H^2(X',\Z))=H^2(X',\Z)^G$ (\cite[Definition 6.1.2]{Weibel})
we have $$E_2^{0,2}=H^0(G,H^2(X',\Z))=H^2(X',\Z)^G,\
E_2^{3,0}=H^3(G,H^0(X',\Z))=H^3(G,\Z),\mbox{ and so}$$
\begin{eqnarray*}\label{formula: Einfty02 second}
E_{\infty}^{0,2}=\ker(H^2(X',\Z)^G\ra H^3(G,\Z)).
\end{eqnarray*}
\textbf{Computation of $E^2$.} By the Theorem \ref{theorem:
convergence of E2pq to Hp+q(S/G,A)} applied to $S=X'$ we obtain
that $E_2^{p,q}\Rightarrow H^2(X'/G,\Z)$, so $E^2=H^2(X'/G,\Z)=H^2(Y',\Z)$.\\
\textbf{Computation of $E_1^2$.} To compute $E_1^2$ we consider
the exact sequence $0\ra E_{2}^{2}\ra E_1^2\ra E_{\infty}^{1,1}\ra
0$ from the diagram \eqref{formula: diagrammone}. By
\eqref{formula: Ep1=0} $H^1(X',\Z)=0$ and so
$E_{2}^{1,1}=0$ but then $E_{\infty}^{1,1}=0$.\\
Since $E_{\infty}^{2,0}=E_2^2/E_3^2$ and $E_3^2=0$, $E_2^2\simeq
E_{\infty}^{2,0}$. The sequence becomes $0\ra E^{2,0}_{\infty}\ra
E_1^2\ra 0$ and then $E_1^2\simeq
E_{\infty}^{2,0}$.\\
Since the spectral sequence is a first quadrant spectral sequence
$E_3^{2,0}=E_{\infty}^{2,0}$ ($E_r^{2,0}=E_{\infty}^{2,0}$ if
$r>\max\{p,q+1\})$. Moreover
$$E_3^{2,0}=\frac{\ker(d_2^{2,0}:E_2^{2,0}\ra E_2^{4,-1})}{\mbox{Im}(d_2^{0,1}:E_2^{0,1}\ra E_2^{2,0})}=E_2^{2,0},$$
because $E_2^{0,1}=0$ by \eqref{formula: Ep1=0}. So
$$E_1^2\simeq E_{\infty}^{2,0}\simeq E_2^{2,0}=H^2(G,\Z).$$

The horizontal sequence of \eqref{formula: diagrammone} becomes
\begin{eqnarray}\label{formula: THE SEQUENCE}0\ra H^2(G,\Z)\ra H^2(Y',\Z)\stackrel{\theta^*}{\ra} \ker(H^2(X',\Z)^G\ra
H^3(G,\Z))\ra 0.\end{eqnarray} The exact sequence \eqref{formula:
THE SEQUENCE} implies that the following sequence is exact:
\begin{eqnarray*}
0\ra H^2(G,\Z)\ra H^2(Y',\Z)\stackrel{\theta^*}{\ra}
H^2(X',\Z)^G.\end{eqnarray*}
\begin{rem}\label{rem: when the map os surjective}{\rm If $H^3(G,\Z)=0$ then the sequence
$$0\ra H^2(G,\Z)\ra H^2(Y',\Z)\stackrel{\theta^*}{\ra}  H^2(X',\Z)^G\ra 0.$$
is exact and in particular the map $$\theta^*:H^2(Y',\Z)\ra
H^2(X',\Z)^G$$
is onto.\\
If $G=\Z/n\Z$ is cyclic, then $H^3(G,\Z)$ is given by
\cite[Example 6.2.3]{Weibel}:
\begin{eqnarray*}
H^n(\Z/n\Z,\Z)=\left\{\begin{array}{ll} \Z&\mbox{ if }
n=0\\
0&\mbox{ if } n=2k+1,\ k\in\N\\
\Z/n\Z&\mbox{ if } n=2k, k\in\N.
\end{array}\right.
\end{eqnarray*}
Hence the map $H^2(Y',\Z)\ra H^2(X',\Z)^G$ is surjective and the
exact sequence is:
$$0\ra \Z/n\Z\ra H^2(Y',\Z)\ra H^2(X',\Z)^G\ra 0.$$
This proves the assertion of Nikulin that $\theta^*$ is surjective
in case $G$ is a cyclic group.\erem}\end{rem}

More in general we want to compute the sequence \eqref{formula:
THE SEQUENCE} for all the finite abelian groups acting
symplectically on a K3 surface (the complete list is given in
\cite[Proposition 1.7]{Nikulin symplectic}).

\begin{prop}\label{prop: sequences for all the group G}
For each finite abelian group $G$ acting sympletically on $X$ we
have the following exact sequences:
$$0\ra G\ra H^2(Y',\Z)\ra ker\left(H^2(X',\Z)^G\stackrel{d_3^{2,0}}{\ra}V_G\right),
\mbox{ where }$$
$$\begin{array}{ll}
V_G=\left\{\begin{array}{lll}0&\mbox{ if
}&G=\Z/n\Z\\
\Z/2\Z&\mbox{ if }&G=(\Z/2\Z)^2\\
(\Z/2\Z)^3&\mbox{ if }&G=(\Z/2\Z)^3\\
(\Z/2\Z)^6&\mbox{ if }&G=(\Z/2\Z)^4
\end{array}\right.&
V_G=\left\{\begin{array}{lll}\Z/3\Z&\mbox{ if
}&G=(\Z/3\Z)^2\\
\Z/4\Z&\mbox{ if }&G=\Z/2\Z\times\Z/4\Z\\
\Z/4\Z&\mbox{ if }&G=(\Z/4\Z)^2\\
\Z/6\Z&\mbox{ if }&G=\Z/2\Z\times\Z/6\Z.
\end{array}\right.\end{array}$$
\end{prop}
\bprf The proof consists in the computation of the groups
$H^2(G,\Z)$ and $H^3(G,\Z)$. This computation is based on the
formulas (cf.\ \cite[Exercise 6.1.8]{Weibel})
\begin{eqnarray*}
H^n(G_1\times G_2)\simeq
\left(\bigoplus_{p+q=n}\!\!H^p(G_1,\Z)\otimes
H^q(G_2,\Z)\right)\oplus\left(\bigoplus_{p+q=n+1}\!\!\!
Tor(H^p(G_1,\Z),H^q(G_2,\Z))\right),\end{eqnarray*}  $Tor(\Z,A)=0$
for each abelian group $A$ and $Tor(\Z/n\Z,\Z/m\Z)=\Z/(n,m)\Z$.
\erem
\begin{rem}{\rm Comparing the sequence \eqref{formula: THE SEQUENCE} with the sequences in
\eqref{diagram: sequences on H2X' and H2Y' } one obtains that
$M_G/\oplus_j\Z M_j\simeq G$ (cf.\ also \cite[Theorem 6.3]{Nikulin
symplectic}, \cite[Lemma 2]{Xiao}).\erem}
\end{rem}

\section{Symplectic automorphisms on $Km(A)$}\label{section: group of automorphisms of order two.} Here we
will construct an example of a K3 surface such that the group
$(\Z/2\Z)^4$ acts symplectically on it. This K3 surface is a
Kummer surface, the symplectic automorphisms are induced by
translation by points of order two on the Abelian
surface.\\
We will compute the lattices $\Lambda_{K3}^G$ for the subgroups
$G\subseteq (\Z/2\Z)^4$. In the case $G=\Z/2\Z$ we obtain the same
result as Nikulin, in fact we have proved in Remark \ref{rem: when
the map os surjective} that the map $\theta^*$ is surjective if
$G$ is cyclic. For the other cases, the sublattice
$\theta^*(P_G^{\vee})$ of $\Lambda_{K3}^G$ is {\em not} equal to
$\Lambda_{K3}^G$. Hence $\theta^*$ is not surjective (cf.\
\eqref{formula: theta* is surjective iff}). As a consequence we
find that the discriminant computed by Nikulin is not correct.
\subsection{Preliminaries on
 Kummer surfaces.}
We recall some properties of Kummer surfaces and of Kummer
lattices. 
\begin{defi}\label{defi: Kummer surface} Let $A$ be an Abelian surface. Let
$\iota:A\ra A$ be the involution $\iota:a\mapsto -a$, $a\in A$.
The quotient $A/\iota$ is a surface with sixteen singularities of
type $A_1$ (the image of the sixteen 2-torsion points of $A$). The
desingularization of this quotient is the K3 surface $Km(A)$,
called \textbf{Kummer surface of $A$}.\\
The minimal primitive sublattice of $H^2(Km(A),\Z)$ containing the
sixteen rational curves resolving the singularities of $A/\iota$
is called \textbf{Kummer lattice}.
\end{defi}
\begin{prop}{\rm \cite[Appendix to section 5, Lemma 4]{Pjateckii Safarevic torelli theorem K3}} The Kummer lattice is a rank sixteen negative
definite lattice. Its discriminant lattice is isomorphic to the
one of $U(2)^{\oplus 3}$ (where $U(2)$ is the lattice obtained
from $U$ by multiplying the bilinear form by two).\end{prop}

Let $A[2]$ be the set of the 2-torsion points on $A$. We fix an
isomorphism $A[2]\simeq (\Z/2\Z)^4$ and we write $p_a\in A[2]$ for
the point corresponding to $a\in(\Z/2\Z)^4$. 
The surface $\widetilde{A}$ is the blow up of $A$ in the 2-torsion
points, so there are sixteen exceptional curves $E_{a}$ on
$\widetilde{A}$ which correspond to the points $p_{a}$. These
sixteen curves are sent in the sixteen rational curves $K_{a}$ of
the Kummer lattice by the map $\pi: \widetilde{A}\ra
\widetilde{A}/\widetilde{\iota}$, where $\widetilde{\iota}$ is the
involution induced on $\widetilde{A}$ by the involution $\iota$ on
$A$ described in Definition \ref{defi: Kummer surface}. We have
the following commutative diagram:
\begin{eqnarray}\label{formula: diagramma Abelian and Kummer}\begin{array}{rrclr}
\{p_{a}\}=A[2]\subset&A&\stackrel{\gamma}{\la}&\widetilde{A}&\supset \{E_{a}\}\\
&\downarrow&\circlearrowleft&\downarrow\pi&\\
Sing(A/\iota)\subset&A/\iota&\la&Km(A)&\supset\{K_{a}\}
\end{array}\end{eqnarray}

We can thus associate a 2-torsion point $p_{a}$ on an Abelian
surface to the corresponding curve $K_{a}$ in the Kummer lattice
of its Kummer surface.
We use the following convention:\\ $\bullet$ if $W$ is an affine
subspace of the affine space $A[2]$, $\bar{K}_W$ is the class
$\frac{1}{2}\sum_{a\in W}K_{a}$ and $\hat{K}=\frac{1}{2}\sum_{a\in
(\Z/2\Z)^4}K_{a}=\bar{K}_{A[2]}$;\\
$\bullet$ $W_i=\{a=(a_1,a_2,a_3,a_4)\in A[2]\mbox{ such that
}a_i=0\}$,
$i=1,2,3,4$;\\
$\bullet$ $W_{i,j}=\{a=(a_1,a_2,a_3,a_4)\in A[2]\mbox{ such that
}a_i+a_j=0\}$,
$1\leq i<j\leq4$;\\
$\bullet$ $V_{i,j}=\{0,\alpha_i,\alpha_j,\alpha_i+\alpha_j\}$,
where $\alpha_1=(1,0,0,0)$, $\alpha_2=(0,1,0,0)$,
$\alpha_3=(0,0,1,0)$, $\alpha_4=(0,0,0,1)$ and $1\leq i<j\leq4$.
\begin{prop}{\rm (\cite[Chapter VIII, Section 5]{bpv}, \cite[Appendix 5]{Pjateckii Safarevic torelli theorem K3})}\label{prop: generator of Kummer lattice}
A set of generators, over $\Z$, of the Kummer lattice is made up
of the sixteen classes: $\hat{K}$, $\bar{K}_{W_1}$,
$\bar{K}_{W_2}$, $\bar{K}_{W_3}$, $\bar{K}_{W_4}$, $K_{a}$, with
$a=(0,0,0,0)$, $(0,0,1,1)$, $(0,1,0,1)$, $(1,0,0,1)$, $(0,1,1,0)$,
$(1,0,1,0)$, $(1,1,0,0)$, $(1,0,0,0)$, $(0,1,0,0)$, $(0,0,1,0)$,
$(0,0,0,1)$.
\end{prop}
\begin{defi}\cite[Definition 5.3]{morrison}\label{latticenik}
The {\bf Nikulin lattice} is an even lattice $N$ of rank eight
generated by $\{N_i\}_{i=1}^8$ and $\hat N=\frac{1}{2}\sum N_i$,
with bilinear form induced by $N_i\cdot N_j=-2\delta_{ij}.$
\end{defi}
\begin{rem}{\rm The Nikulin lattice is the lattice $M_G$ (cf.\ Definition \ref{defi: MG})
where $G$ is generated by an involution, the $N_i$ correspond to
curves arising from the singularities of the quotient of a K3
surface by a symplectic involution.\erem}
\end{rem}
\begin{rem}{\rm For each $i$ the classes $K_a$, $a\in W_i$, and $\bar{K}_{W_i}$ generate a Nikulin lattice.}\erem
\end{rem}

\begin{prop}{\rm (\cite[Chapter VIII, Section 5]{bpv}, \cite[Appendix 5]{Pjateckii Safarevic torelli theorem K3})}\label{prop: discriminant Kummer lattice}
Let $A_K=K^{\vee}/K$ be the discriminant group of the lattice $K$
and $q_{A_K}:A_K\ra \Q/2\Z$ be the quadratic form on $A_K$ induced
by the pairing on $K$.\\The discriminant group $A_K$ is made up of
63 elements and the zero $0_{A_K}$: 35 of them are of type
$\bar{K}_{V}$, for linear subspaces $V\simeq (\Z/2\Z)^2\subset
A[2]$, and $q_{A_{K}}(\bar{K}_V)=0$, the other 28 are of type
$\bar{K}_{V+V'}$, for linear subspaces $V$ and $V'$, such that
$V\cap V'=\{0\}$, and $q_{A_K}(\bar{K}_{V+V'})=1$. There are three
orbits of $O(q_{A_K})$ on $A_K$: $\{0\}$, $\{v\in A_K:
q_{A_K}(v)\equiv 0 \mod 2\Z\}$ and $\{v\in A_K: q_{A_K}(v)\equiv 1
\mod 2\Z\}$.
\end{prop}
\begin{rem}\label{rem: generators of dicriminant of K}{\rm The
classes $\bar{K}_{V_{1,2}}$, $\bar{K}_{V_{3,4}}$,
$\bar{K}_{V_{1,3}}$, $\bar{K}_{V_{2,4}}$, $\bar{K}_{V_{1,4}}$,
$\bar{K}_{V_{2,3}}$ generate the discriminant group $A_K$.\\
Moreover if $e_i$, $f_i$, $i=1,2,3$ are the standard basis of
$U(2)^{\oplus 3}$, then $\bar{K}_{V_{1,2}}+e_1/2$,
$\bar{K}_{V_{1,3}}+e_2/2$, $\bar{K}_{V_{1,4}}+e_3/2$,
$\bar{K}_{V_{3,4}}+f_1/2$, $\bar{K}_{V_{2,4}}+f_2/2$,
$\bar{K}_{V_{2,3}}+f_3/2$ together with the classes generating the
Kummer lattice (cf.\ Proposition \ref{prop: generator of Kummer
lattice}) are 22 vectors in $(K\oplus U(2)^{\oplus
3})\otimes_{\Z}\Q$ which generate a unimodular even lattice with
signature $(3,19)$, so a lattice isometric to
$\Lambda_{K3}$.}\erem
\end{rem}
\begin{defi} Let $x_i$, $i=1,2,3,4$, be the real coordinates of an Abelian surface
$A=(\R/\Z)^4\simeq \C^2/\Lambda$. Let
$\omega_{i,j}:=\pi_*(\gamma^*(dx_i\wedge dx_j))$, with the
notation of the diagram \eqref{formula: diagramma Abelian and
Kummer}.
\end{defi}
\begin{rem}{\rm The classes $\omega_{i,j}$, $1\leq i<j\leq 4$ generate $U(2)^{\oplus
3}\simeq K^{\perp}\subset H^2(Km(A),\Z)$.\\
\binf since $H^4(X,\Z)\simeq \Z$, with $dx_1\wedge dx_2\wedge
dx_3\wedge dx_4\mapsto 1$, and $$(dx_i\wedge dx_j)\wedge
(dx_h\wedge dx_k)\mapsto\left\{\begin{array}{l} \pm 1 \mbox{ if
}\{i,j,h,k\}=\{1,2,3,4\}\\
0 \mbox{ otherwise,}\end{array}\right.$$ the six forms $dx_1\wedge
dx_2$, $dx_3\wedge dx_4$, $dx_1\wedge dx_3$, $dx_2\wedge dx_4$,
$dx_1\wedge dx_4$, $dx_2\wedge dx_3$ generate three copies of the
lattice $U$ respectively. Hence they form a basis for
$H^2(A,\Z)\simeq U^{\oplus 3}$.\\
The forms $dx_i\wedge dx_j$ on $A$ are invariant under $\iota^*$,
and so they generate $H^2(A,\Z)^{\iota^*}=H^2(A,\Z)$. Let us
consider the lattice $H^2(A,\Z)^{\iota^*}$ as sublattice of
$H^2(\widetilde{A},\Z)$ (it is exactly the sublattice orthogonal
to the classes of the exceptional curves $E_a$). Since
$\pi_*(H^2(A,\Z))^{\iota^*}\simeq H^2(A,\Z)^{\iota^*}(2)$ (cf.\
\cite[Lemma 3.1]{morrison}), the lattice generated by
$\omega_{i,j}=\pi_*(\gamma^*(dx_i\wedge dx_j))$, is isometric to
$U(2)^{\oplus 3}$ and it is a sublattice of $H^2(Km(A),\Z)$.
Moreover the classes $\omega_{i,j}$, are orthogonal to the curves
of the Kummer lattice in $H^2(Km(A),\Z)$, and so they generate
$K^{\perp}\subset H^2(Km(A),\Z)$.}\erem
\end{rem}
\begin{rem}\label{rem: divisible calsses and omega}{\rm We fix a complex
structure on $A$ such that the complex coordinates are
$z_1=x_1+ix_2$ and $z_2=x_3+ix_4$, and we choose the lattice such
that $\C^2/\Lambda\simeq \C/\Gamma\times \C/\Gamma'$ where
$\Gamma$ and $\Gamma'$ are lattices in $\C$. This defines two
elliptic curves $C:=\C/\Gamma$ and $C':=\C/\Gamma'$ (with real
coordinates $(x_1,x_2)$ and $(x_3,x_4)$ respectively) and hence
two classes $[C]$ and $[C']$ in $H_2(A,\Z)$. By Poincar\'e duality
$H_2(A,\Z)\simeq H^2(A,\Z)^*\simeq H^2(A,\Z)$, so $[D]$
corresponds to a two form $\mu$ such that $\int_D\nu=
\int_S\mu\wedge \nu$ for all $\nu\in H^2(S,\R)$. As $$\int_C
dx_i\wedge dx_j=\left\{\begin{array}{rl}1&\mbox{ if
}i=1, j=2,\\
-1&\mbox{ if
}i=2, j=1,\\
0&\mbox{ otherwise}\end{array}\right.$$ and $\int_{A} dx_h\wedge
dx_k\wedge dx_1\wedge dx_2= 1$ if and only if $h=3$ and $k=4$, the
class $[C]$ corresponds to $dx_3\wedge dx_4$. Similarly
$[C']=dx_1\wedge dx_2$. On each of these curves there are four
2-torsion points. With the notation of \eqref{formula: diagramma
Abelian and Kummer}, $\gamma^*(C)=\widetilde{C}+\sum_{a\in C[2]}
E_a$ where $\widetilde{C}$ is the strict transform of the $C$ and
$E_a$ are the exceptional curves over the four 2-torsion points of
$A[2]$ on the curve $C$ (i.e. $a=(a_1,a_2,0,0)$). By Remark
\ref{rem: generators of dicriminant of K} we know that
$(u/2+\bar{K}_{V_{1,2}})\in\Lambda_{K3}$ for a certain $u\in
U(2)^{\oplus 3}$ with $u^2=0$. The curves $K_{a_1,a_2,0,0}$
correspond to the four 2-torsion points on the curve $C$, and the
class $[C]$ is $dx_3\wedge dx_4$, so $u=\omega_{3,4}$. This
implies that $\omega_{3,4}/2+\bar{K}_{V_{1,2}}$ is a class in
$H^2(Km(A),\Z)$. So we can restate Remark \ref{rem: generators of
dicriminant of K} in the following, more precise, way: a set of
generators of the lattice $H^2(Km(A),\Z)$ is given by the sixteen
classes generating the Kummer lattice and by the six classes
$\omega_{1,2}/2+\bar{K}_{V_{3,4}}$,
$\omega_{3,4}/2+\bar{K}_{V_{1,2}}$,
$\omega_{1,3}/2+\bar{K}_{V_{2,4}}$,
$\omega_{2,4}/2+\bar{K}_{V_{1,3}}$,
$\omega_{1,4}/2+\bar{K}_{V_{2,3}}$,
$\omega_{2,3}/2+\bar{K}_{V_{1,4}}$. \erem}\end{rem}

\subsection{The surfaces $Km(A)$ and $Km(A/\langle b\rangle )$} $ $\\
Let $A$ be an Abelian surface. Let $b$ be a point of order two on
$A$ and $\langle b\rangle=\{0,b\}$ be the group generated by $b$.
\begin{rem}\label{rem: two torsion points on A quotient b}{\rm The surface $A/\langle b\rangle$ is again an Abelian surface.
The 2-torsion points on $A/\langle b\rangle$ are the images of the
points $r$ on $A$ such that $2r\in\langle b\rangle$. We observe
that the points $r$ such that $2r=b$ are 4-torsion points on $A$.
Hence the sixteen points in $(A/\langle b\rangle)[2]$ are the
eight points in the image of the 2-torsion points of $A$ and the
eight points in the image of the 4-torsion points $r$ such that
$2r=b$.\erem}\end{rem}

The translation by the point $b$ of order two on $A$ induces an
involution on $Km(A)$ and on $\widetilde{Km(A)}$. Since the
translation by the point $b$ acts as the identity on the
holomorphic two forms on $A$ and since the holomorphic two form on
$Km(A)$ is induced by the one on $A$, the involution induced on
$Km(A)$ by the translation by $b$ fixes the holomorphic two form
on $Km(A)$ and hence is symplectic. We have the following
commutative diagram:
$$
\begin{array}{rcclccr}
&&A&-\!\!\!-\!\!\!-\!\!\!-\!\!\!-\!\!\!-\!\!\!-\!\!\!-\!\!\!-\!\!\!-\!\!\!-\!\!\!-\!\!\!-\!\!\!-\!\!\!-\!\!\!-\!\!\!-\!\!\!-\!\!\!-\!\!\!-\!\!\!-\!\!\!-\!\!\!-\!\!\!-\!\!\!\rightarrow&A/\langle b\rangle&&\\
&\swarrow&&\searrow\ \ \ \ \ \ \ \ \ \ \ \ \ \ \ \ \ \ \ \ \ \ \ \ \ \ \ \ \ \ \ \ \ \swarrow&&\searrow&\\
A/\iota&&&\ \ \ \widetilde{A}\ \ \ \ \ \ \ \ \ \ \ \ \ \ \ \ \ \ \
(A/\langle b\rangle )/\iota
&&&\widetilde{A/\langle b\rangle }\\
&\searrow&&\swarrow\ \ \ \ \ \ \ \ \ \ \ \ \ \ \ \ \ \ \ \ \ \ \ \ \ \ \ \ \ \ \ \ \ \searrow&&\swarrow&\\
&&Km(A)&\leftarrow\!\!\!-\!\!\!-\!\!\!-
\widetilde{Km(A)}-\!\!\!-\!\!\!-\!\!\!-\!\!\!-\!\!\!-\!\!\!-\!\!\!-\!\!\!-\!\!\!-\!\!\!-\!\!\!-\!\!\!\rightarrow&Km(A/\langle
b\rangle )=\widetilde{Km(A)/\langle b\rangle }&&
\end{array}
$$
where $\widetilde{Km(A)}$ is the blow up of $Km(A)$ in the points
of $Km(A)$ with a non-trivial stabilizer for the action of the
involution induced by the translation by $b$ and
$\widetilde{Km(A)/\langle b\rangle }$ is the smooth quotient of
$\widetilde{Km(A)}$ by the action this involution on $Km(A)$.

\begin{rem}\label{rem: branch locus of cover between Kummer}{\rm
Let $x\in A$ such that its image $\overline{x}\in A/\iota$ is a
point with a non trivial stabilizer (i.e. a fixed point) with
respect to the translation by $b$. Then we have
$\overline{x+b}=\overline{x}$ and so $x+b=\pm x$. Since $b$ is of
order two, $2x=b$ and so $x$ has order four. There are 16 such
points $x$ on $A$ and so on $A/\iota$ we find $16/2=8$ fixed
points. After blow up of the images of the 2-torsion points in
$A/\iota$ we obtain $Km(A)$, in particular there are 8 fixed
points in $Km(A)$ (with respect to the involution induced on
$Km(A)$ by the translation on $A$). Note that none of the curves
of the Kummer lattice on $Km(A)$ are fixed. Thus the branch locus
of $\widetilde{Km(A)}\ra Km(A/\langle b\rangle)$ are the eight
curves on $Km(A/\langle b\rangle)$ corresponding to the eight
points in $(A/\langle b \rangle)[2]$ which are images of $x$ with
$2x=b$.\erem}\end{rem}


\subsection{The group $G=(\Z/2\Z)^4$.}$ $\\
Let $G=(\Z/2\Z)^4$ be the group of symplectic automorphisms on
$Km(A)$ induced by the group of the translation by points of order
two on $A$. As $A/A[2]\simeq A$ the desingularization of the
quotient of $Km(A)$ by $G$ is again $Km(A)$. The image of the
point $p_{a}$ on $A$ under the map $\varphi:A\ra A/ A[2]$ is the
point $p_{2a}$, so in particular the image of $p_{a}$ with
$a\in(\Z/2\Z)^4=A[2]$ under the map $A\ra
A/ A[2]$ is the point $p_{0}$.\\

\textbf{The lattice $M_{(\Z/2\Z)^4}$.} Since the 2-torsion points
$p_{a}$ are sent to $p_0$, the curves $K_{a}$ on $Km(A)$ are sent
to the curve $K_0$ on $\widetilde{Km(A)/G}=Km(A)$. The branch
locus of the cover $\widetilde{Km(A)}\ra\widetilde{Km(A)}/G$ can
be found as in Remark \ref{rem: branch locus of cover between
Kummer}. In fact the stabilizer, with respect to the action of
$G$, of a point $\overline{x}\in A/\iota$ is either trivial or the
group $\{0,b\}$, (and then $2x=b$) for a certain $b\in
A[2]\backslash\{0\}$. Since there are fifteen 2-torsion points,
there are fifteen different stabilizer groups. For each of these
stabilizer groups there are eight points with that group as
stabilizer (cf.\ Remark \ref{rem: branch locus of cover between
Kummer}). Hence there are $8\cdot 15=120$ points with a non
trivial stabilizer on $A/\iota$. The quotient by the group $G$
identifies the eight points with the same stabilizer group, hence
the branch locus of the cover $\widetilde{Km(A)}\ra
\widetilde{Km(A)/G}$ is made up of the fifteen curves of the
Kummer lattice $K_a$, $a\neq 0$.
By definition $M_{(\Z/2\Z)^4}$ is the minimal primitive sublattice
of $H^2(Y,\Z)$ containing these fifteen curves,
\begin{eqnarray*}\langle K_{a}, a\in
(\Z/2\Z)^4\backslash \{0\}\rangle\stackrel{f.\
i.}{\hookrightarrow} M_{(\Z/2\Z)^4},\end{eqnarray*} where
the inclusion has a finite index.\\

\textbf{The lattice $P_{(\Z/2\Z)^4}=M_{(\Z/2\Z)^4}^{\perp}$.}
Since $M_{(\Z/2\Z)^4}\subset K$, the lattice
$K^{\perp}=U(2)^{\oplus 3}$ is contained in
$M_{(\Z/2\Z)^4}^{\perp}=P_{(\Z/2\Z)^4}$. Moreover the curve $K_0$
is contained in $M_{(\Z/2\Z)^4}^{\perp}$. So $\Z K_0\oplus
U(2)^{\oplus 3}$ is contained, with finite
index, in $P_{(\Z/2\Z)^4}$.\\
There are no classes which are in $P_{(\Z/2\Z)^4}$ but not in $\Z
K_0\oplus U(2)^{\oplus 3}$, in fact there are no classes in the
dual of the lattice $K\oplus U(2)^{\oplus 3}$ involving only one
curve of the Kummer lattice (cf.\ Remark \ref{rem: generators of
dicriminant of K} and Remark \ref{rem: divisible calsses and
omega}), so
\begin{eqnarray*}
\begin{array}{cl}
P_{(\Z/2\Z)^4}&=\left\langle K_0,
\omega_{1,2},\omega_{3,4},\omega_{1,3},\omega_{2,4},\omega_{1,4},\omega_{2,3}\right\rangle\\
&\simeq\langle -2\rangle\oplus U(2)\oplus U(2)\oplus U(2).
\end{array}
\end{eqnarray*}
\textbf{The lattice $\theta^*(P_{(\Z/2\Z)^4}^{\vee})$.} The dual
lattice of $P_{(\Z/2\Z)^4}$ is made up of the classes $c\in
P_{(\Z/2\Z)^4}\otimes\Q$ such that the intersection product
$c\cdot r\in \Z$ for each $r\in P_{(\Z/2\Z)^4}$. Then
\begin{eqnarray*}
P_{(\Z/2\Z)^4}^{\vee}=\left\langle K_0/2,
\omega_{1,2}/2,\omega_{3,4}/2,\omega_{1,3}/2,\omega_{2,4}/2,\omega_{1,4}/2,\omega_{2,3}/2\right\rangle.
\end{eqnarray*}
To determine $\theta^*(P_G^{\vee})$, note that
$\theta^*(K_0)=\sum_{a\in(\Z/2\Z)^4}K_{a}$. The map $\vartheta$,
and the hence also map $\theta$ (cf.\ \eqref{formula: vartheta}),
is induced by the map $A\stackrel{\cdot 2}{\ra} A$, so its action
on $\omega_{i,j}$ is induced by the action of $\cdot 2$ on
$dx_i\wedge dx_j$. The map $\cdot 2:dx_i\wedge dx_j\ra
d(2x_i)\wedge d(2x_j)=4dx_i\wedge dx_j$, and so
$\theta^*(\omega_{i,j})=4\omega_{i,j}$. We observe that
$(\theta^*(\omega_{i,j}),\theta^*(\omega_{h,k}))=(4\omega_{i,j},4\omega_{h,k})=|G|(\omega_{i,j},\omega_{h,k})$
as predicted.\\
By the descriptions of $P_{(\Z/2\Z)^4}^{\vee}$ and $\theta^*$ we
obtain:
\begin{eqnarray*}
\begin{array}{c}\theta^*(P_{(\Z/2\Z)^4}^{\vee})=\left\langle
\hat{K},
2\omega_{1,2},2\omega_{3,4},2\omega_{1,3},2\omega_{2,4},2\omega_{1,4},2\omega_{2,3}\right\rangle\\\simeq\langle-8\rangle\oplus
U(8)\oplus U(8)\oplus U(8).
\end{array}
\end{eqnarray*}
\textbf{The lattice $\Lambda_{K3}^{(\Z/2\Z)^4}$.} The group $G$
acts trivially on the $\omega_{i,j}$ and so the $\omega_{i,j}$ are
contained in $\Lambda_{K3}^{(\Z/2\Z)^4}$. Moreover each class
$K_a$ of the Kummer lattice is sent into $K_{a+a'}$ by $a'\in G$,
so no class $K_a$ is fixed, but the sum of all these classes is
fixed. So $\hat{K}$ is contained in $\Lambda_{K3}^{(\Z/2\Z)^4}$.
The classes $\omega_{i,j}$ generate the lattice $U(2)^{\oplus 3}$
which is primitive in $\Lambda_{K3}$, and so
\begin{eqnarray*}\label{invariant Z2Z^4}
\begin{array}{cl}
\Lambda_{K3}^{(\Z/2\Z)^4}&=\left\langle \hat{K},
\omega_{1,2},\omega_{3,4},\omega_{1,3},\omega_{2,4},\omega_{1,4},\omega_{2,3}\right\rangle\\
&\simeq\langle-8\rangle\oplus U(2)\oplus U(2)\oplus U(2).
\end{array}
\end{eqnarray*}
Comparing the discriminant or observing the classes $\omega_{i,j}$
it is clear that $\theta^*(P_{(\Z/2\Z)^4}^{\vee})$ is not
isomorphic to $\Lambda_{K3}^{(\Z/2\Z)^4}$ but it is contained with
a finite index in $\Lambda_{K3}^{(\Z/2\Z)^4}$,
$$\theta^*(P_{(\Z/2\Z)^4}^{\vee})\stackrel{f.i}{\hookrightarrow} \Lambda_{K3}^{(\Z/2\Z)^4}\ \ \ \ \ \mbox{  and  }\ \ \ \ \ [\Lambda_{K3}^{(\Z/2\Z)^4}:\theta^*(P_{(\Z/2\Z)^4}^{\vee})]=2^6.$$
In particular $\theta^*$ is not surjective in this case.
\subsection{The group $G=(\Z/2\Z)^2$}$ $\\
Let us consider the subgroup $G_{2,2}$ of $(\Z/2\Z)^4$ induced by
the translation by two 2-torsion points. The quotient $A/G_{2,2}$
is an Abelian surface. With $A=(\R/\Z)^4$, we choose the map $A\ra
A/G_{2,2}$ to correspond to multiplication by two on the first two
coordinates, $x_1$ and $x_2$, on $A$. We observe that $A/G_{2,2}$
is not isomorphic to $A$ in general.\\
The subgroup $G_{2,2}$
induces the group $G=(\Z/2\Z)^2$ of symplectic automorphisms on
$Km(A)$.\\

\textbf{The lattice $M_{(\Z/2\Z)^2}$.} It is the minimal primitive
sublattice of $H^2(Y,\Z)$ containing the twelve curves $K_{a},
a=(a_1,a_2,a_3,a_4)\in (\Z/2\Z)^4,
(a_1,a_2)\in(\Z/2\Z)^2\backslash\{0\}$ which are the branch locus
of the
cover $\widetilde{Km(A)}\ra \widetilde{Km(A)/G}.$\\

\textbf{The lattice $P_{(\Z/2\Z)^2}=M_{(\Z/2\Z)^2}^{\perp}$.} By
the computation of $M_{(\Z/2\Z)^2}$, the $16-12=4$ remaining
curves of $K$ are in $P_{(\Z/2\Z)^2}$, in fact:
$$\langle
\omega_{1,2},\omega_{3,4},\omega_{1,3},\omega_{2,4},\omega_{1,4},\omega_{2,3},K_0,
K_{0,0,0,1}, K_{0,0,1,0},
K_{0,0,1,1}\rangle\stackrel{f.i.}{\hookrightarrow}
P_{(\Z/2\Z)^2},$$ with a finite index. Fixing a complex structure
on $A$ as in Remark \ref{rem: divisible calsses and omega}, we
know that $\bar{K}_{V_{3,4}}+\omega_{1,2}/2\in H^2(Y,\Z)$. Then:
$$P_{(\Z/2\Z)^2}=\left\langle \bar{K}_{V_{3,4}}+\omega_{1,2}/2,\omega_{3,4},\omega_{1,3},\omega_{2,4},\omega_{1,4},\omega_{2,3},K_0, K_{0,0,0,1}, K_{0,0,1,0}, K_{0,0,1,1}\right\rangle.$$
\textbf{The lattice $\theta^*(P_{(\Z/2\Z)^2}^{\vee})$.} The dual
lattice of $P_{(\Z/2\Z)^2}$ is
\begin{eqnarray*}
\begin{array}{c}
P_{(\Z/2\Z)^2}^{\vee}=\left\langle
\omega_{1,2}/2,\omega_{1,3}/2,\omega_{1,4}/2,\omega_{2,3}/2,\omega_{2,4}/2,\left(\omega_{3,4}+K_0\right)/2\right.,\\
\left.\left(K_0+K_{0,0,0,1}\right)/2,\left(K_0+K_{0,0,1,0}\right)/2,
\left(K_0+K_{0,0,1,1}\right)/2,K_0
\right\rangle,\mbox{ and}\\
\\
\theta^*(P_{(\Z/2\Z)^2}^{\vee})= \left\langle
2\omega_{1,2},\omega_{1,3},\omega_{1,4},\omega_{2,3},\omega_{2,4},\left(\omega_{3,4}/2+\bar{K}_{V_{1,2}}\right),\bar{K}_{W_3},
\bar{K}_{W_4}, \bar{K}_{W_{3,4}}, 2\bar{K}_{V_{1,2}}\right\rangle.
\end{array}
\end{eqnarray*}

\textbf{The lattice $\Lambda_{K3}^{(\Z/2\Z)^2}$.} As before the
classes $\omega_{i,j}$ are in $\Lambda_{K3}^{(\Z/2\Z)^2}$. The
group $G$ acts as the identity on the third and fourth coordinates
of $A$, so for $(a_3,a_4)\in(\Z/2\Z)^2$ classes
$K_{0,0,a_3,a_4}+K_{0,1,a_3,a_4}+K_{1,0,a_3,a_4}+K_{1,1,a_3,a_4}$
are in $\Lambda_{K3}^{(\Z/2\Z)^2}$. Over $\Q$,
$\Lambda_{K3}^{(\Z/2\Z)^2}$ is generated by $\omega_{i,j}$ and
$\sum_{(a_1,a_2)\in(\Z/2\Z)^2} K_{a_1,a_2,a_3,a_4}$, $(a_3,a_4)\in
(\Z/2\Z)^2$. The sum of eight curves $K_{a}$ which correspond to a
hyperplane are divisible by two in the Kummer lattice (cf.\
Proposition \ref{prop: generator of Kummer lattice}), so in
particular they are elements of $\Lambda_{K3}$. Moreover
$\bar{K}_{V_{1,2}}+\omega_{3,4}/2$ is contained in $\Lambda_{K3}$
(cf.\ Remark \ref{rem: divisible calsses and omega}). Then the
primitive lattice in $\Lambda_{K3}$ generated by $\omega_{i,j}$
and $\sum_{(a_1,a_2)\in(\Z/2\Z)^2} K_{a_1,a_2,a_3,a_4}$,
$a_3,a_4\in (\Z/2\Z)^2$ is
\begin{eqnarray*}\label{equation: invariant 2,2}
\begin{array}{c}
\Lambda_{K3}^{(\Z/2\Z)^2}=\left\langle\omega_{1,2},\omega_{1,3},
\omega_{2,4},\omega_{1,4},\omega_{2,3},\left(\omega_{3,4}/2+\bar{K}_{V_{1,2}}\right),
\bar{K}_{W_3}, \bar{K}_{W_4}, \bar{K}_{W_{3,4}},2\bar{K}_{V_{1,2}}
\right\rangle.\end{array}\end{eqnarray*} Comparing this lattice
with $\theta^*(P_{(\Z/2\Z)^2} ^{\vee})$ it is clear that they are
not isomorphic and
$$\theta^*(P_{(\Z/2\Z)^2}^{\vee})\stackrel{f.i}{\hookrightarrow} \Lambda_{K3}^{(\Z/2\Z)^2}\ \ \ \ \ \mbox{  and  }\ \ \ \ \ [\Lambda_{K3}^{(\Z/2\Z)^2}:\theta^*(P_{(\Z/2\Z)^2}^{\vee})]=2.$$
\subsection{$G=(\Z/2\Z)^3$}$ $\\
Similar computations can be done in case $G=(\Z/2\Z)^3$. Let
$G_{2,2,2}$ be the group of translations by three independent
points of order two on $A$, inducing the group $G$  of symplectic
automorphisms on $Km(A)$. We may assume that the quotient map
$A\ra A/G_{2,2,2}$ corresponds to a multiplication by two on the
first three real coordinates of
$A$.\\
In this case one obtains:
\begin{eqnarray*}\label{invariant Z2Z^3}
\begin{array}{ccl}
P_G=\left\langle
\omega_{1,2},\omega_{3,4},\omega_{1,3},\omega_{2,4},\omega_{1,4},\omega_{2,3},K_0,
K_{0,0,0,1}\right\rangle,\\
\\
\begin{array}{cl}
\theta^*(P_{(\Z/2\Z)^3}^{\vee})&= \left\langle\right.
2\omega_{1,2},\omega_{3,4},2\omega_{1,3},\omega_{2,4},\omega_{1,4},2\omega_{2,3},
\bar{K}_{W_4}, \bar{K}_{(\Z/2\Z)^4-W_4}\left.\right\rangle\\
&\simeq U(4)\oplus U(4)\oplus U(4)\oplus \langle
-4\rangle\oplus\langle-4\rangle,
\end{array}\\
\\
\begin{array}{cl}
\Lambda_{K3}^{(\Z/2\Z)^3}&=\left\langle\right.\omega_{1,2},\omega_{3,4},\omega_{1,3},
\omega_{2,4},\omega_{1,4},\omega_{2,3},\bar{K}_{W_4},
\bar{K}_{(\Z/2\Z)^4-W_4}\left.\right\rangle\\
&\simeq U(2)\oplus U(2)\oplus U(2)\oplus \langle
-4\rangle\oplus\langle-4\rangle.
\end{array}
\end{array}
\end{eqnarray*}
Comparing $\Lambda_{K3}^{(\Z/2\Z)^3}$ and
$\theta^*(P_{(\Z/2\Z)^3}^{\vee})$, one obtains
$[\Lambda_{K3}^{(\Z/2\Z)^3}:\theta^*(P_{(\Z/2\Z)^3}^{\vee})]=2^3.$\\

\subsection{$G=\Z/2\Z$}$ $\\
In this case $G$ is generated by a symplectic involution (which is
induced by the translation by a point of order two). In Remark
\ref{rem: when the map os surjective} we proved that in this case
the map $H^2(Y',\Z)\ra H^2(X,\Z)^G$ is surjective, and then
$\Lambda_{K3}^{\Z/2\Z}=\theta^*(P_{\Z/2\Z}^{\vee})$. This result
can also be obtained with the method used above. In fact let $G_2$
be the translation by a point of order two on $A$ such that the
map $A\ra A/G_2$ corresponds to a multiplication by two on the
first real coordinate of $A$. Let $G$ be the group generated by
the involution induced by $G_2$ on $Km(A)$. Then
$$
\begin{array}{l}
\begin{array}{lcl}
P_{\Z/2\Z}&=&\langle \bar{K}_{V_{3,4}}+\omega_{1,2}/2,
\omega_{3,4}, \bar{K}_{V_{2,4}}+\omega_{1,3}/2,\omega_{2,4},
\bar{K}_{V_{2,3}}+\omega_{1,4}/2,\omega_{2,3}, \bar{K}_{W_1},\\&&\
\ \ \
K_{0,0,0,1},K_{0,0,1,0},K_{0,1,0,0},K_{0,0,1,1},K_{0,1,0,1},K_{0,1,1,0},K_{0,1,1,1}\rangle,\\
\\
\theta^*(P_{\Z/2\Z}^{\vee})&= &\langle
\omega_{1,2},\omega_{3,4},\omega_{1,3},\omega_{2,4},\omega_{1,4},\omega_{2,3},
\bar{K}_{W_2},\bar{K}_{W_3},\bar{K}_{W_4},\hat{K},
K_{0,0,1,1}+K_{1,0,1,1},\\&&\ \ \ \ K_{0,1,0,1}+K_{1,1,0,1},K_{0,1,1,0}+K_{1,1,1,0},K_{0,1,1,1}+K_{1,1,1,1}\rangle,\\
\\
\Lambda_{K3}^{\Z/2\Z}&=&\langle
\omega_{1,2},\omega_{3,4},\omega_{1,3},\omega_{2,4},\omega_{1,4},\omega_{2,3},
\bar{K}_{W_2},\bar{K}_{W_3},\bar{K}_{W_4},\hat{K},
K_{0,0,1,1}+K_{1,0,1,1},\\&&\ \ \ \
K_{0,1,0,1}+K_{1,1,0,1},K_{0,1,1,0}+K_{1,1,1,0},K_{0,1,1,1}+K_{1,1,1,1}\rangle.
\end{array}
\end{array}
$$
The last two lattices are the same, as proved before.
\begin{rem} {\rm We observe that for each group $G=(\Z/2\Z)^i$,
$i=2,3,4$ the lattice $\theta^*(P^{\vee}_G)$ has a finite index in
$\Lambda_{K3}^G$ and $\Lambda_{K3}^G/\theta^*(P^{\vee}_G)\simeq
V_G$, where $V_G$ is the group described in Proposition \ref{prop:
sequences for all the group G}.\erem}\end{rem}
\begin{rem}{\rm In \cite{symplectic not prime} the discriminant of
the lattice $\Lambda_{K3}^G$ is computed for each finite abelian
group acting symplectically on a K3 surface. For each non cyclic
group it is different from the one computed by Nikulin. Comparing
the two discriminant we always obtain that the index of the
inclusion of the lattice computed by Nikulin in our lattice is
exactly the order of $V_G$, with $V_G$ as in Proposition
\ref{prop: sequences for all the group G}.\erem}\end{rem}

\end{document}